\documentclass[11pt,notitlepage]{article}
\usepackage{amsmath}
\usepackage{enumitem}
\usepackage{amssymb}
\usepackage{enumitem}
\usepackage{fullpage}
\usepackage{microtype}
\usepackage{url}
\usepackage[small,sf]{caption}
\usepackage[font=sf]{floatrow}
\usepackage{verbatim}
\usepackage{comment}
\usepackage{graphicx}\usepackage{amsfonts}\usepackage{amsthm}\usepackage{lscape}
\usepackage{booktabs}
\usepackage{amscd}
\usepackage[usenames,dvipsnames]{color}
\usepackage[colorlinks,linkcolor=blue]{hyperref}
\usepackage{chngcntr}
\counterwithin{table}{section}
\hypersetup{citecolor=blue}

\usepackage{listings}
\usepackage{color}

\definecolor{dkgreen}{rgb}{0,0.6,0}
\definecolor{gray}{rgb}{0.5,0.5,0.5}
\definecolor{mauve}{rgb}{0.58,0,0.82}

\allowdisplaybreaks[1]

\makeatletter
\renewcommand{\l@section}{\@dottedtocline{1}{1.5em}{2.3em}}
\makeatother

\theoremstyle{remark}

\newcommand{\sfcaption}[1]{\caption{{\sf #1}}}

\newcommand{\M}{\mathcal M}

\newcommand{\cc}{\mathbb C}
\newcommand{\chat}{\widehat \cc}

\newcommand{\rr}{\mathbb R}

\newcommand{\zz}{\mathbb Z}

\begin{document}
\author{Ronen E. Mukamel\footnote{The research for this paper was supported in part by grant DMS-1103654 from the National Science Foundation.}}
\title{Visualizing the unit ball for the Teichm\"uller metric}
\date{October 9, 2014}
\maketitle
\begin{abstract}
We describe a method to compute the norm on the cotangent space to the moduli space of Riemann surfaces associated to the Finsler Teichm\"uller metric.  Our method involves computing the periods of abelian double covers and is easy to implement for Riemann surfaces presented as algebraic curves using existing tools for approximating period matrices of plane algebraic curves.  We illustrate our method by depicting the unit sphere in the cotangent space to moduli space at a particular surface of genus zero with five punctures and by corroborating the proof of a theorem of Royden's for our example.
\end{abstract}
\begin{figure}
  \begin{center}
   \includegraphics[scale=0.6]{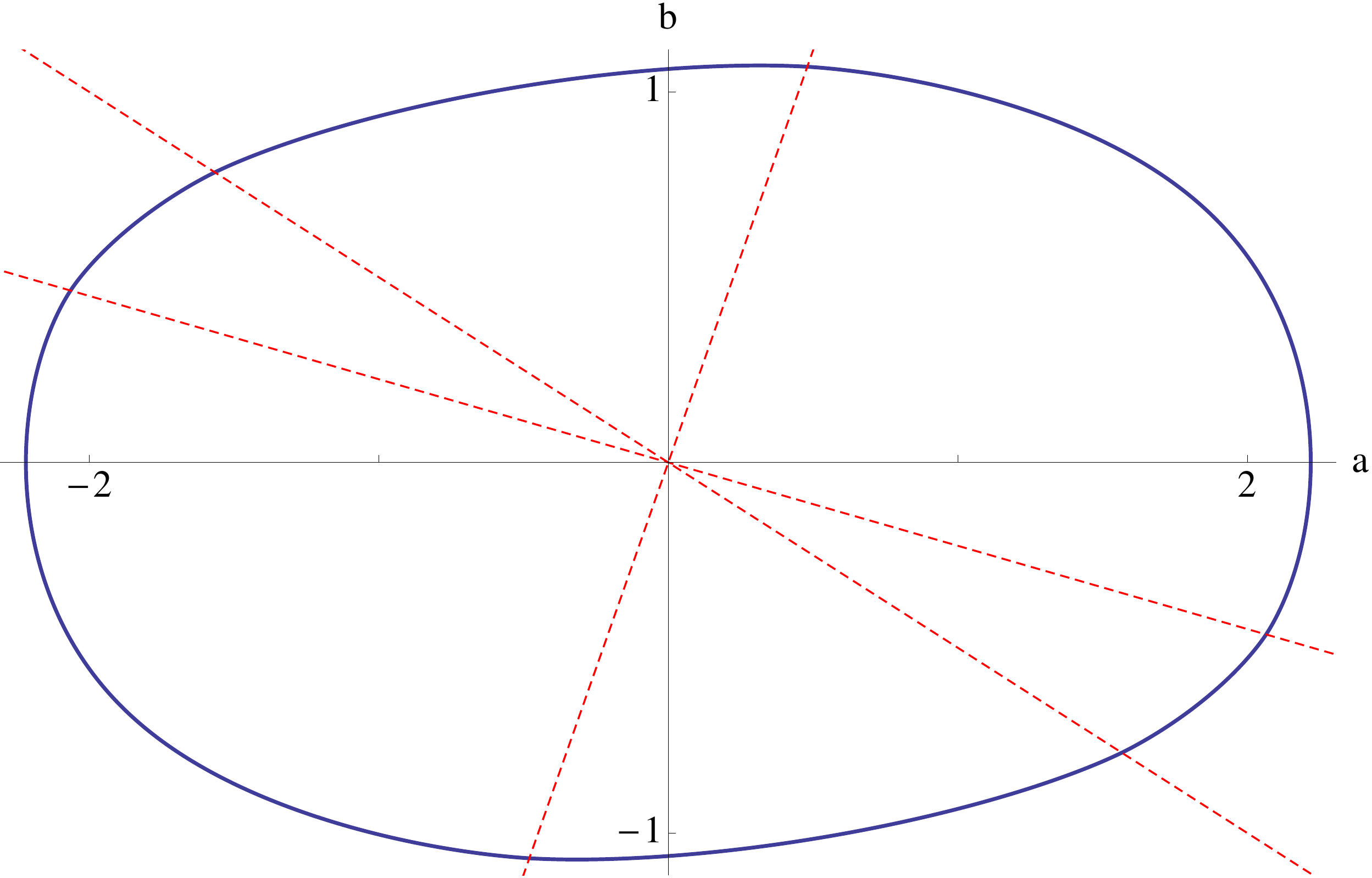}
\end{center}
\sfcaption{\label{fig:unitball} 
Using our algorithm, we sample the unit sphere $U$ (solid) defined by Equations \ref{eqn:exampleh} and \ref{eqn:unitreal} in the vector space $Q(X)$ for a particular surface $X \in \M_{0,5}$.  The sphere $U$ is smooth except where it meets the differentials in $Q(X)$ with four poles and no zeros (dashed).}
\end{figure}
\begin{figure}
  \begin{center}
   \includegraphics[scale=0.45]{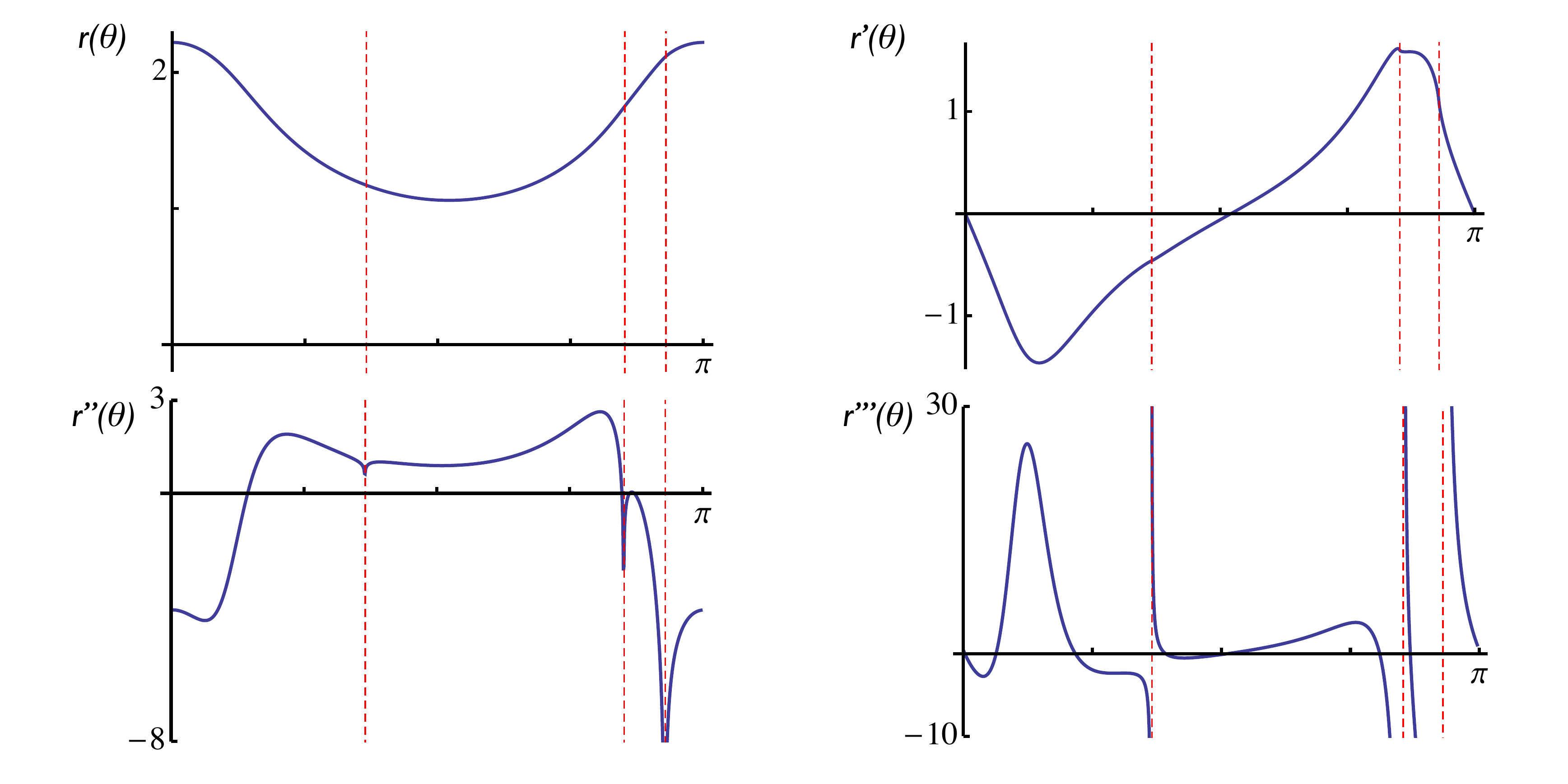}
\end{center}
\sfcaption{\label{fig:derivatives} The sphere $U$ in Figure \ref{fig:unitball} in polar coordinates (solid, top-left) is smooth except at angles corresponding to forms with no zeros (dashed) where $r''(\theta)$ is unbounded (cf. \cite[Par. 2.4]{earlekra:holomorphic} and discussion below).}
\end{figure}

The purpose of this paper is to describe a method to compute the norm on the cotangent space to the moduli space of Riemann surfaces whose study was pioneered by Royden in \cite{royden:automorphisms}.  Royden's norm is important in Teichm\"uller theory because its dual gives rise to a Finsler metric on moduli space which is equal to the Teichm\"uller metric defined using quasiconformal mappings.  Royden's norm and the Teichm\"uller metric are connected with many areas of mathematics and their study has celebrated applications to the classification of mapping classes \cite{bers:extremal,thurston:surfacediffeos}, the dynamics of rational maps \cite{douadyhubbard:rationalfunctions}, the complex geometry of moduli space \cite{royden:automorphisms,hubbard:sections} and polygonal billiards \cite{veech:ngon,kerckhoff:ergodicity}.  Our method to compute Royden's norm at a particular Riemann surface $X$, summarized in Figure \ref{fig:algorithm}, involves computing the periods of certain double covers of $X$.  Existing tools for approximating period matrices of plane algebraic curves yield simple (numerical) implementations of our algorithm for Riemann surfaces presented as algebraic curves.  In Figure \ref{fig:genuszeroimplementation} we give a short implementation for Riemann surfaces of genus zero.  We demonstrate our method by depicting the unit sphere in the cotangent space to moduli space at a particular Riemann surface of genus zero with five punctures (Figure \ref{fig:unitball}) and by giving striking corroboration of a theorem of Royden's for our example (Figure \ref{fig:derivatives}).

\paragraph{Quadratic differentials and singular flat metrics.}  Let $\M_{g,n}$ denote the moduli space of Riemann surfaces of finite type $(g,n)$.  A point in $\M_{g,n}$ is a Riemann surface $X$ which is isomorphic to the complement of $n$ points on a closed Riemann surface $\overline{X}$ of genus $g$.  A {\em holomorphic quadratic differential} $q$ on $X$ is a holomorphic section of the symmetric square of the cotangent bundle on $X$.  At each $x \in X$, such a differential restricts to a quadratic form $q_x : T_xX \to \cc$ and, when $q_x \neq 0$, the absolute value of $q_x$ defines a norm $|q_x| : T_xX \to \rr$.  These norms give rise to an area form $|q|$ on $X$ associated to a flat metric with cone singularities at the zeros of $q$.

\paragraph{Integrable quadratic differentials.}  For a quadratic differential $q$ on $X$, we define $\| q \|$ to be the area of the associated singular flat metric, i.e.
\begin{equation}
\label{eqn:rnorm}
\| q \| = \int_X |q|.
\end{equation}
We call $q$ {\em integrable} if $\| q\|$ is finite.  Equivalently, $q$ is integrable if it extends to a meromorphic differential on $\overline{X}$ whose poles are simple.  

The vector space $Q(X)$ is the vector space of all integrable holomorphic quadratic differentials on $X$ and Royden's norm is the $L^1$-norm on $Q(X)$ defined by Equation \ref{eqn:rnorm}.  When $m = 3g-3+n$ is positive, $Q(X)$ has complex dimension $m$ by Riemann-Roch and $Q(X)$ is naturally isomorphic to the cotangent space $T^*_X\M_{g,n}$ by Teichm\"uller theory.  Royden showed that the dual of the norm $\| \cdot \|$ gives rise to a Finsler metric on $\M_{g,n}$ which equals the Teichm\"uller metric defined via quasiconformal geometry and the Kobayashi metric $\M_{g,n}$ inherits as a complex orbifold \cite{royden:automorphisms}.

\paragraph{Polygonal presentations of quadratic differentials.}  An important way to present an integrable quadratic differential $q$ is to specify a finite collection of polygons in $\cc$ and glue parallel sides of equal length together by maps of the form $z \mapsto \pm z+c$.  The resulting quotient $\overline{X}$ inherits the structure of a closed Riemann surface from $\cc$ and a meromorphic quadratic differential $q$ from the form $dz^2$ which is automatically integrable. In fact, the area of $|q|$ is the sum of the areas of the polygons used to construct $q$ and therefore the norm $\| q \|$ is elementary to compute.  A challenge for understanding the norm $\| \cdot \|$ from this perspective is that, although every integrable holomorphic quadratic differential admits a polygonal presentation, it is difficult to determine whether two polygonally presented differentials are differentials on isomorphic Riemann surfaces.  

\paragraph{Algebraic quadratic differentials.}  Alternatively, one can present $X$ as an algebraic curve and $q$ as an algebraic differential.  Fix three polynomials $f,g,h \in \cc[x,y]$ with $f$ irreducible and $h$ coprime to $f$.  These polynomials determine a closed Riemann surface $\overline{X}$ and a meromorphic quadratic differential $q$ on $\overline{X}$ via the equations
\begin{equation}
\label{eqn:algdiffl}
 \cc\!\left(\overline{X}\right) \cong \cc(x)[y] / (f) \mbox{ and } q = \frac{g(x,y)}{h(x,y)} dx^2.
\end{equation}
One can use resultants to locate the poles of $q$ and Puiseux series to give conditions on $f$, $g$ and $h$ ensuring that $q$ is integrable.  

For a particular $X \in \M_{g,n}$ presented as an algebraic curve, it is often possible to describe the entire vector space $Q(X)$ as algebraic differentials.  For instance, if $\overline{X}$ is the Riemann sphere (e.g. $f(x,y) = y$), $h \in \cc[x]$ is a polynomial with simple roots and $X \subset \overline{X}$ is the complement of $Z(h)$ the set of zeros of $h$, then the vector space $Q(X)$ satisfies
\begin{equation}
\label{eqn:QXS}
Q(X) = \left\{ \frac{g(x)}{h(x)} dx^2 : g \in \cc[x] \mbox{ with } \deg(g) \leq \deg(h) - 4 \right\}.
\end{equation}
For curves of higher genus, the vector space $Q(X)$ consists of differentials of the form in Equation \ref{eqn:algdiffl} where $f$ and $h$ are fixed and the coefficients of $g$ satisfy certain linear conditions.

\paragraph{Abelian square roots.}  Now let $\Omega\!\left(\overline{X}\right)$ be the space of abelian differentials (i.e. holomorphic one-forms) on $\overline{X}$.  A form $\omega \in \Omega\!\left(\overline{X}\right)$ restricts to a linear function $\omega_x : T_x\overline{X} \to \cc$ for each $x \in \overline{X}$ and the square of $\omega_x$ is a quadratic form on $T_x\overline{X}$.  These quadratic forms vary holomorphically in $x$ and, in this way, the square of $\omega$ can be viewed as a holomorphic quadratic differential $q = \omega^2$ on $\overline{X}$.  Forms which are squares can be polygonally presented using transition maps of the restricted form $z \mapsto z+c$.  The algebraic differential defined by Equation \ref{eqn:algdiffl} is the square of an abelian differential if and only if it is integrable and $g/h$ is a square in the field $\cc(x)[y]/(f)$.  

\paragraph{Norms of squares of abelian differentials.}  When $q = \omega^2$ happens to be the square of an abelian differential $\omega$, the norm $\| q \|$ is easily computed from the periods of $\omega$.  For any symplectic basis $\left< a_1, \dots, a_g, b_1, \dots, b_g \right>$ of $H_1\left(\overline{X},\zz\right)$, the norm of $q$ satisfies
\begin{equation}
\label{eqn:normfromperiods}
 \| q \| = \sum_{j=1}^g \operatorname{Im}\left( \overline{x}_j y_j \right) \mbox{ where } x_j = \int_{a_j} \omega \mbox{ and } y_j = \int_{b_j} \omega.
\end{equation}
Equation \ref{eqn:normfromperiods} follows from the identity $|q| = (i/2) \cdot \omega \wedge \overline{\omega}$ and Stokes' theorem.

\paragraph{Abelian double covers.}  For an arbitrary quadratic differential $q \in Q(X)$, we define $\overline{X}_q$ to be the normalization of the closure of 
\begin{equation}
\label{eqn:dblcover}
\left\{ (x,\omega_x) : x \in X \mbox{ and } \omega_x^2 = q_x \right\} \subset T^* \overline{X}.
\end{equation}
The projection $\pi : \overline{X}_q \to \overline{X}$ is degree two and branched at the zeros and poles of $q$ of odd order.  The pullback $\pi^*q$ is the square of an abelian differential $\omega_q$ on $\overline{X}_q$ and satisfies $\| \pi^* q \| = 2 \| q \|$.  We call the form $\omega_q$ the {\em abelian double cover of $q$}.  

When $X$ and $q$ are given by Equation \ref{eqn:algdiffl}, the surface $\overline{X}_q$ and abelian double cover $\omega_q$ are defined by the equations
\begin{equation}
\label{eqn:algcandblcover}
\cc\!\left( \overline{X}_q \right) \cong \cc\!\left( \overline{X} \right)[z] / (z^2 - g \cdot h) \mbox{ and } \omega_q \cong \frac{g}{z} \cdot dx.
\end{equation}
Note that, when $q$ is a square, $\overline{X}_q$ has two irreducible components.

\begin{figure}
\setlength{\fboxsep}{0.25in}
\begin{center}
\fbox{
\parbox{5in}{
{\bf Algorithm for computing Royden's norm}
\newline
\newline \noindent {\bf Input:} An integrable, algebraic quadratic differential $q$
\newline \noindent {\bf Output:} The norm $\| q \|$
\newline
\begin{enumerate}[leftmargin=*,label=\textbf{(\arabic*)}]
\item Form the surface $\overline{X}_q$ and the abelian double cover $\omega_q$ as in Equation \ref{eqn:algdiffl},
\item Compute the periods of $\omega_q$ (e.g. using \cite{deconinckvhoejj:Riemannmatrices}),
\item Compute $\|\omega_q^2 \|$ using Equation \ref{eqn:normfromperiods}, and
\item Return $\| q \|= \frac{1}{2} \|\omega_q^2 \|$.
\end{enumerate}
}
}
\end{center}
\sfcaption{\label{fig:algorithm} Royden's norm for a particular quadratic differential can be computed from the periods of the abelian double cover.}
\end{figure}
\paragraph{Algorithm.}  From the fact that every quadratic differential is double covered by the square of an abelian differential and Equation \ref{eqn:normfromperiods}, computing Royden's norm reduces to computing the periods of abelian differentials.  An algorithm for computing such periods for differentials presented algebraically is described in \cite{deconinckvhoejj:Riemannmatrices}.  To compute the periods of $\omega_q$ defined by Equation \ref{eqn:algcandblcover}, one might:
\begin{enumerate}
\item compute the critical values $B$ of the map $x : \overline{X}_q \to \chat$,
\item choose a graph $G \subset \chat \setminus B$ such that inclusion is a homotopy equivalence, and
\item numerically integrate $\omega_q$ along edges of $\widetilde G = x^{-1}(G)$.
\end{enumerate}
The numerical integration in (3) is easiest for arcs that avoid $x^{-1}(B)$, so a natural choice for $G$ is the Voronoi diagram of $\chat$ relative to $B$.  Because the inclusion $\widetilde G \to \overline{X}_q$ induces a surjection on first homology, the integrals computed in (3) determine the periods of $\omega_q$.  Variants of this algorithm have been implemented in {\sf MAGMA} for hyperelliptic curves (cf. \cite{vwamelen:analyticjacobians}) and in {\sf Maple} for arbitrary plane algebraic curves.

Our algorithm for computing Royden's norm is now clear and is summarized in Figure \ref{fig:algorithm}.  Similar ideas can likely be used to create a robust tool to give polygonal presentations of algebraic quadratic differentials, and such a tool would be very useful.

\paragraph{Implementation in genus zero.}  Now suppose $h \in \cc[x]$ is a polynomial with simple roots and $X$ is the complement of the zeros $Z(h)$ of $h$ in the Riemann sphere.  The vector space $Q(X)$ is given by Equation \ref{eqn:QXS}.  For any particular $q = (g/h) \cdot dx^2$ in $Q(X)$, the surface $\overline{X}_q$ is the hyperelliptic curve defined by the equation
\begin{equation}
\label{eqn:weierstrass}
z^2 = g(x) \cdot h(x)
\end{equation}
and the one-form
\begin{equation}
\label{eqn:hyperellipticsqrt}
\omega_q = g(x) \cdot \frac{dx}{z}
\end{equation}
is an abelian square root of $\pi^*q$.  The tools in {\sf MAGMA} related to analytic Jacobians compute period matrices of hyperelliptic curves, yielding the short implementation of our algorithm for genus zero $X$ in Figure \ref{fig:genuszeroimplementation}.  

\paragraph{Genus zero example.}  Now specialize further to the (arbitrarily chosen) case
\begin{equation}
\label{eqn:exampleh}
 h(x) = x^5-2 x^4-12 x^3-8 x^2+52 x +24 \mbox{ and } Z(h) = \left\{ 2, 2 \pm \sqrt{6}, -2 \pm \sqrt{-2} \right\}.
\end{equation}
As in the previous paragraph, let $X \in \M_{0,5}$ be the complement of $Z(h)$ in the Riemann sphere.  The vector space $Q(X)$ is two dimensional and, in Figure \ref{fig:unitball}, we use our algorithm to sample the (real) unit sphere 
\begin{equation}
\label{eqn:unitreal} U = \left\{ (a,b) \in \rr^2 : \left\| \frac{a+bx}{h(x)} \cdot dx^2 \right\|= 1 \right\} \subset Q(X).
\end{equation}
Figure \ref{fig:unitball} is generated using 1000 samples each computed using 100 digits of precision.  On the machine (3.4 GHz, 8GB RAM) used to generate Figure \ref{fig:unitball}, individual samples compute instantaneously and the entire sampling process takes approximately two minutes.

\begin{figure}
\lstset{tabsize=4}
\begin{lstlisting}
function QDNorm(qD)
//Input: A list qD = [g,h] where g and h are polynomials with coefficients in a complex field; g and h should satisfy deg(g) <= deg(h) - 4, GCD(g,h) = 1 and Discriminant(h) neq 0.
//Output: the norm of g/h*dx^2
	A := AnalyticJacobian(qD[1]*qD[2]); 
	dimA := Dimension(A); degg := Degree(qD[1]); 
	pers := &+[Coefficient(qD[1],j-1)*BigPeriodMatrix(A)[j] : j in [1..(degg+1)]];
	return &+[Im(pers[j]*Conjugate(pers[j+dimA])) : j in [1..dimA]]/2;
end function;
\end{lstlisting}
\sfcaption{\label{fig:genuszeroimplementation} Implementation of our algorithm in {\sf MAGMA} (version V2.20-3) for genus zero quadratic differentials.}
\end{figure}

\paragraph{Smoothness of Royden's norm.}   In \cite{royden:automorphisms}, Royden proved several amazing facts about the norm $\| \cdot \|$ on $Q(X)$ for $X \in \M_{g,0}$ and the associated Finsler metric $d_T$ on $\M_{g,0}$.  The analogous facts for $\M_{g,n}$ are established by similar means in \cite{earlekra:holomorphic,earlekra:isometries}.  These authors show that the metric $d_T$ equals both the Teichm\"uller metric on $\M_{g,n}$ arising from quasiconformal geometry and the Kobayashi metric $\M_{g,n}$ inherits as a complex orbifold.  They also show that the biholomorphic automorphisms of the Teichm\"uller space universal covering $\M_{g,n}$ is the mapping class group when $m = 3g-3+n > 1$.

A key ingredient in Royden's argument is to show that, when $\dim_\cc(Q(X)) > 1$, the normed vector space $\left( Q(X),\| \cdot \|\right)$ determines $X$ up to isomorphism.  This fact is established by showing that the norm on $Q(X)$, when restricted to lines passing through $q \in Q(X)$, is of class $C^1$ and has a H\"older continuous derivative whose exponent is determined by the orders of the zeros of $q$.  For other results related to the smoothness of Royden's norm and consequences for the geometry of moduli space see e.g. \cite{earle:c1,rees:c2,rees:c2e,antonakoudis:tdisks}.

We illustrate this phenomenon for our example in Figure \ref{fig:derivatives} by plotting the unit sphere $U$ defined by Equation \ref{eqn:unitreal} in polar coordinates along with the first few derivatives of the Euclidean radius $r(\theta)$ as a function of angle along $U$.  Derivatives were computed using successive difference quotients.  The third derivative of $r(\theta)$ clearly acquires singularities at the three angles corresponding to forms with four simple poles and no zeros corresponding to the three real roots of $h(x)$.  In fact, according to \cite[Par. 2.4]{earlekra:holomorphic}, the norm $\| \cdot \|$ is of class $C^{1,\alpha}$ for every $\alpha < 1$.  The second derivative of $\| \cdot \|$ is unbounded at the differentials with no zeros, although its growth is slow near such differentials.

\paragraph{Acknowledgments.}  The author would like to thank C. T. McMullen, A. Epstein and the referee for helpful suggestions.  

\bibliography{numericalqds}

\bigskip
\noindent {\sc Ronen E. Mukamel \\
Department of Mathematics \\ 
Rice University, MS 136 \\
6100 Main St. \\
Houston, TX 77005}

\end{document}